\documentclass[12pt]{amsart}

\usepackage{amssymb}
\usepackage{epsf}

\newtheorem{thm}{Theorem}[section]
\newtheorem{prop}[thm]{Proposition}
\newtheorem{lem}[thm]{Lemma}
\newtheorem{cor}[thm]{Corollary}
\newtheorem{definition}[thm]{Definition}
\newtheorem{prob}[thm]{Problem}
\newtheorem{rem}[thm]{Remark}
\newenvironment{Remark}{\begin{rem}\sl}{\end{rem}}

\newcommand{\abs}[1]{|#1|}
\def\silentfootnote#1{{\let\thefootnote\relax\footnotetext{#1}}}

\def\u#1{{\bf u}_{#1}}
\def\hu#1{\hat{{\bf u}}_{#1}}

\def\hx{\hat{\bf x}}
\def\hy{\hat{\bf y}}
\def\hz{\hat{\bf z}}
\def\Sym{{\mathcal S}}
\def\Q{{\mathbb Q}}
\def\C{{\mathbb C}}
\def\N{{\mathbb N}}
\def\Z{{\mathbb Z}}
\def\Mon{{\mathcal M}}
\def\Non{{\mathcal N}}
\def\a{\alpha}
\def\b{\beta}
\def\c{\gamma}
\def\d{\delta}

\begin{document}

\title[A Monoid for the universal $k$-Bruhat order]{A Monoid for 
the universal $k$-Bruhat order} 

\author{Nantel Bergeron \and Frank Sottile}

\address{Department of Mathematics and Statistics\\
        York University\\
        North York, Ontario M3J 1P3\\
	CANADA}
\email[Nantel Bergeron]{bergeron@mathstat.yorku.ca}
\address{Department of Mathematics\\
        University of Toronto\\
        100 St.~George Street\\
	Toronto, Ontario  M5S 3G3\\
	CANADA}
\email[Frank Sottile]{sottile@math.toronto.edu}
\date{2 July 1997}
\thanks{Bergeron supported in part by CRM, MSRI, and NSERC}
\thanks{Sottile supported in part by NSF grant DMS-9022140,
NSERC grant  OGP0170279, and CRM}
\subjclass{05E15, 14M15, 05E05}
\keywords{Bruhat order, nil-Coxeter monoid, flag manifold, 
Grassmannian}

\begin{abstract}
Structure constants for the multiplication of Schubert polynomials by Schur
symmetric polynomials are known to be related to the enumeration of chains
in a new partial order on $\Sym_\infty$, the universal $k$-Bruhat
order. 
Here we present a monoid $\Mon$ for this order. 
We show that $\Mon$ is analogous to the nil-Coxeter monoid for the weak
order on $\Sym_\infty$. For this, we develop the theory of reduced
sequences for $\Mon$. 
We use these sequences to give a combinatorial description of the structure
constants above. 
We also give a combinatorial proof of some of the symmetry relations
satisfied by these constants.
\end{abstract}

\maketitle

\section{Introduction}

Let $\Sym_\infty$ denote the infinite symmetric group consisting
of permutations of $\{1,2,\ldots\}$ which fix all but finitely many numbers.
In their approach to the Schubert calculus for flag manifolds, 
Lascoux and Sch\"utzenberger~\cite{LS82a,LS82b,LS83,LS89}
defined Schubert polynomials 
${\mathfrak S}_u\in \Z[x_1,x_2,\ldots]$, a homogeneous basis indexed by 
permutations $u\in\Sym_\infty$.
By construction, the degree of ${\mathfrak S}_u$ is the length, 
$\ell(u)$, of $u$. 
We refer the reader to \cite{Macdonald91} for an interesting detailed 
account of Schubert polynomials and double Schubert polynomials. 
This construction has been extended to quantum Schubert polynomials for the
manifolds of complete flags~\cite{FGP,KM} and for manifolds of partial
flags~\cite{C-F}. 
In~\cite{Fulton97}, W. Fulton generalizes all of these constructions.

It is a famous open problem to understand the multiplicative
structure constants for the Schubert polynomials and any of their
generalizations. 
This would provide an understanding of some Gromov-Witten invariants.
From algebraic geometry, the structure constants $c_{uv}^w$ defined by the
identity 
  $$
{\mathfrak S}_u{\mathfrak S}_v\ =\ 
\sum_{w\in \Sym_\infty} c_{uv}^w {\mathfrak S}_w
$$
are known to be positive integers, and in some cases they reduce to the 
Littlewood-Richardson coefficients. 
A general combinatorial construction or bijective formula for the $c_{uv}^w$
is not known.

It is believed that $c_{uv}^w$ counts the number of chains from $u$ to $w$ 
in the Bruhat order which satisfy conditions imposed by $v$ \cite{BS97a}. 
In particular, if $v$ is a Grassmannian permutation with descent in $k$,
then one can restrict the chains to a suborder: 
the $k$-Bruhat order $\le_k$ on 
$\Sym_\infty$~\cite{LS83,Sottile96,BS97a}.
In \cite{BS97a}, a study of $\le_k$ leads to a new partial order 
$\preceq$ on $\Sym_\infty$ which we
call the {\sl universal
$k$-Bruhat} order. This order is ranked and has the property that
a nonempty interval $[u,w]_k$ in a $k$-Bruhat order is
isomorphic to the interval $[1,wu^{-1}]_\preceq$ in the universal order
(independent of $k$).
Every interval in Young's lattice is an interval in this universal order.
The first aim of this paper is to present a monoid $\Mon$ that
describes the chain structure of the universal $k$-Bruhat order.

The monoid $\Mon$ has a $0$ and generators $\u{\a\b}$ indexed by 
integers $0<\a<\b$, subject
to the relations
$$	
\begin{array}{clrclll}
(1)&&\u{\b\c}\u{\c\d}\u{\a\c}&\equiv&\u{\b\d}\u{\a\b}\u{\b\c},\hfill&&
        \hbox{if $\a<\b<\c<\d$},\hfill\\
(2)\hfill&& \hfill \u{\a\c}\u{\c\d}\u{\b\c}&\equiv&
\u{\b\c}\u{\a\b}\u{\b\d},\hfill&&
        \hbox{if $\a<\b<\c<\d$},\hfill\\
(3)\hfill&&  \hfill \u{\a\b}\u{\c\d}&\equiv&\u{\c\d}\u{\a\b}, \hfill&&
        \hbox{if $\b<\c$ or $\a<\c<\d<\b$},\hfill\\
(4)\hfill&&  \hfill   \u{\a\c}\u{\b\d}&\equiv& \u{\b\d}\u{\a\c}\ 
\equiv \  0,\hfill&&
        \hbox{if $\a\le \b<\c\le\d$},\hfill\\
(5)\hfill&&  \hfill  \u{\b\c}\u{\a\b}\u{\b\c}
&\equiv&\u{\a\b}\u{\b\c}\u{\a\b}\ \equiv\ 0,\hfill&&
        \hbox{if $\a<\b<\c$}.
\end{array} \eqno{(1.1)}
$$	

The relation between $\Mon$ and the order $\preceq$ on $\Sym_\infty$ is
obtained via a faithful 
representation of $\Mon$ as linear operators on $\Q\Sym_\infty$. 
Let $\ell_{\u{}}$ denote the rank function of $\preceq$.
Let $(\a\,\,\,\b)\in\Sym_\infty$ be the transposition that interchanges $\a$
and $\b$.  
We define linear operators $\hu{\a\b}$ by 
$$	
\begin{array}{rcl}
\hu{\a\b}\ \colon\  \Q\Sym_\infty&\longrightarrow& \quad\Q\Sym_\infty,\\
\zeta\quad&\longmapsto&\ \ \rule{0pt}{28pt} \left\{\begin{array}{ll}
(\a\,\,\,\b)\zeta   
&\mbox{ if } \ell_{\u{}}\big((\a\,\,\,\b)\zeta)\big)=\ell_{\u{}}(\zeta)+1,\\
\mbox{ }\\      0& \mbox{ otherwise.}\end{array}\right.
\end{array}  \eqno{(1.2)}
$$	

The main results of Section 3 are summarized in the following theorem.

\begin{thm}\label{Operators} \mbox{ }
\begin{enumerate}
\item[(a)]
The map $\ell_{\u{}}\colon\Sym_\infty\to\N$ is well defined by
$\ell_{\u{}}(\zeta)=\ell(\zeta u)-\ell(u)$ for any $u$ and 
$k$ such that $u\le_k\zeta u$.
\item[(b)] 
The operators $\hu{\a\b}$ satisfy the relations (1.1), and a
composition of operators is characterized by its value at the identity. 
That is
$\hu{\a'_m\b'_m}\cdots\hu{\a'_1\b'_1}=\hu{\a_n\b_n}\cdots\hu{\a_1\b_1}$ 
if and only if
\/$\hu{\a'_m\b'_m}\cdots\hu{\a'_1\b'_1}1=\hu{\a_n\b_n}\cdots\hu{\a_1\b_1}1$.
\item[(c)] 
For ${\bf x}=\u{\a_n\b_n}\cdots\u{\a_2\b_2}\u{\a_1\b_1}\in\Mon$, the map 
${\bf x}\mapsto
\hx=\hu{\a_n\b_n}\cdots\hu{\a_2\b_2}\hu{\a_1\b_1}$
is a faithful representation of $\Mon$.
\item[(d)] 
The following map is a well defined bijection:
\begin{eqnarray*}
\Mon&\longrightarrow& \Sym_\infty\cup\{0\},\\
{\bf x}&\longmapsto& \hx 1.
\end{eqnarray*}
\item[(e)] 
The universal $k$-Bruhat order $\preceq$ on $\Sym_\infty$ is ranked by 
$\ell_{\u{}}$.
We have $\eta\preceq \zeta$ if and only if there exists 
${\bf x}\in \Mon$ such that $\zeta=\hx
\eta$.  The order $\preceq$  satisfies the universal property:
$[u,\zeta u]_k\cong [1,\zeta]_\preceq$ whenever $u\le_k \zeta u$.
In particular $[\eta,\zeta]_\preceq\cong[1,\zeta\eta^{-1}]_\preceq$ 
whenever $\eta\preceq\zeta$.
\item[(f)] 
The set 
$R_{\u{}}(\zeta)=\{\hx:\hx1=\zeta\}$
corresponds to the set of all maximal chains in $[1,\zeta]_\preceq$.
\end{enumerate}
\end{thm}

We call the elements of $R_{\u{}}(\zeta)$ the $\u{}$-reduced sequences of
$\zeta$. 
Parts (a) and (e) of Theorem \ref{Operators} are obtained in \S 3.2 of
\cite{BS97a}. 
In Section 2, we relate Theorem \ref{Operators} to classical results on the
weak order of $\Sym_\infty$ and the nil-Cotexer monoid.

Recall \cite{Macdonald91} that the Schur polynomial 
$S_\lambda(x_1,x_2,\ldots,x_k) ={\mathfrak S}_{v(\lambda,k)}$ for a 
unique Grassmannian permutation $v(\lambda,k)$.
In Theorem E of \cite{BS97a}, we have shown that if 
$c_{uv(\lambda,k)}^w\neq 0$, then $c_{uv(\lambda,k)}^w$
depends only on $\lambda$ and
$\zeta=wu^{-1}$. We can thus define constants $c_\lambda^\zeta$ such that 
$c_{uv(\lambda,k)}^w=c_\lambda^{wu^{-1}}$ whenever $u\le_k w$.
We note that ({\sl cf.} Proposition 1.1 \cite{BS97a})
$$	
	\abs{R_{\u{}}(\zeta)} = \sum_\lambda f^\lambda c_\lambda^\zeta,
	\eqno{(1.3)}
$$	
where $f^\lambda$ is the number of standard Young tableaux of shape $\lambda$.
In Section 4 we give a combinatorial description of the constant 
$c_{uv(\lambda,k)}^w$ using elements of $R_{\u{}}(\zeta)$. 
We use this description to give a combinatorial proof of
many of the symmetry relations given in \cite{BS97a}.
In Section 5 we discuss open problems related to the monoid 
$\Mon$ and the constants
$c_{uv(\lambda,k)}^w$. 

The interested reader may obtain by email from 
{\tt bergerna@mathstat.yorku.ca} or find at {\tt
http://www.math.yorku.ca/Who/Faculty/Bergeron} two appendices. 
One describes a graphical
representation of chains in $\preceq$ which greatly helps visualize the
relations (1.1) and the 
arguments of \S 3. The other describes an insertion correspondence, giving a
bijection between  
${\mathcal H}_{n,1}(\zeta)$ and ${\mathcal H}_{1,n}(\zeta)$. This is related
to one open problem 
described in Section 5.

\section{orders and monoids on $\Sym_\infty$}

The weak order $\le_{weak}$ on $\Sym_\infty$ is the transitive closure of the
following cover relation: for $u,w\in\Sym_\infty$, we say that $w$ covers
$u$ in the weak order if 
$\ell(w)=\ell(u)+1$ and $wu^{-1}$ is a simple transposition $(\a\,\,\a\!+\!1)$.
Maximal chains from the identity to $w\in\Sym_\infty$ correspond to
reduced  sequences for $w$. The nil-Coxeter monoid $\Non$ plays an important
role~\cite{EG,LS82b,Macdonald91} in studying reduced sequences. 
The monoid $\Non$ has a $0$
and generators
$\u{i}$ indexed by integers $i>0$, subject to the
nil-Coxeter relations:
$$	
\begin{array}{rcl}
  \u{\a}\u{\a+1}\u{\a}&\equiv&\u{\a+1}\u{\a}\u{\a+1},\\ 
  \u{\a}\u{\b}&\equiv&\u{\b}\u{\a}, \qquad\qquad\hbox{if $\abs{\a-\b}>1$},\\ 
  \u{\a}\u{\a}&\equiv& 0.
\end{array}	\eqno{(2.1)}
$$	
There is a faithful representation of $\mathcal N$ as linear
operators on the group algebra $\Q\Sym_\infty$. 
For this, let
$$
\begin{array}{rcl}
\hu{\a}\ \colon \ \Q\Sym_\infty&\longrightarrow& \quad\Q\Sym_\infty,\\ 
        \zeta\quad&\longmapsto&\ \ \rule{0pt}{28pt} 
\left\{\begin{array}{lll}
(\a\,\,\,\a\!+\!1)\zeta&\quad&
\mbox{if }\ell\big((\a\,\,\,\a\!+\!1)\zeta\big)=\ell(\zeta)+1,\\  \\ 
0&&\mbox{otherwise.}
\end{array}\right.\end{array}
$$

The following proposition is a reformulation of well known facts about 
reduced sequences of a permutation and the weak order. 
See \cite{Macdonald91} for a proof of most of them.
\vfill \eject

\begin{prop}\label{WeakOp}\mbox{ }
\begin{enumerate}
\item[(a)] The map $\ell\ \colon\ \Sym_\infty\to\N$ is well defined.
\item[(b)] The operators $\hu{\a}$ satisfy the relations~(2.1), and a
composition of operators is characterized by its value at the identity. 
That is
$\hu{\a_n}\cdots\hu{\a_1}=\hu{\b_m}\cdots\hu{\b_1}$ if and only if
$\hu{\a_n}\cdots\hu{\a_1}1=\hu{\b_m}\cdots\hu{\b_1}1$.
\item[(c)] For ${\bf x}=\u{\a_n}\cdots\u{\a_2}\u{\a_1}\in\Non$, the map ${\bf
x}\mapsto
\hx=\hu{\a_n}\cdots\hu{\a_2}\hu{\a_1}$
is a faithful representation of $\Non$.
\item[(d)] The following map is a well defined bijection:
\begin{eqnarray*}
\Non&\longrightarrow& \Sym_\infty\cup\{0\},\\
{\bf x}&\longmapsto& \hx 1.
\end{eqnarray*}
 \item[(e)] The weak order $\le_{weak}$ on $\Sym_\infty$ is ranked by $\ell$.
We have $u\le_{weak} w$ if and only if there exists ${\bf x}\in \Non$ such
that $w=\hx u$. 
Also $[\eta,\zeta]_{weak}\cong[1,\zeta\eta^{-1}]_{weak}$ whenever $\eta
\le_{weak}\zeta$.
\item[(f)] The set 
$R(w)=\{\hx:\hx1=w\}$
corresponds to the set of all maximal chains in $[1,w]_{weak}$.
The elements of $R(w)$ are the reduced sequences of $w$.
\end{enumerate}
\end{prop}

At this point we note the striking resemblance between 
Theorem~\ref{Operators} and Proposition~\ref{WeakOp}. 
The proof of Proposition \ref{WeakOp} relies on the understanding of reduced
sequences. 
For Theorem \ref{Operators}, the order $\preceq$ is new and its chains 
have not been studied previously. 
We develop the elementary theory of the analogue of reduced sequences for
$\preceq$. 

We note that not all orders on $\Sym_\infty$ have such a simple monoid.
In particular, the Bruhat order $\le$ on $\Sym_\infty$ has no known monoid.
Recall that $w$ covers $u$ in the Bruhat order if $\ell(w)=\ell(u)+1$ 
and $wu^{-1}$ is a
transposition
$(\a\,\,\b)$. In fact, very little is known
about the problem of chain enumeration for the Bruhat order.
We believe that a monoid for the Bruhat order would not satisfy 
conditions as simple as
those of Theorem \ref{Operators} and Proposition \ref{WeakOp}. 

The monoid structure for the weak order was
a key factor in the following results.
Under the nil-Coxeter-Knuth relations
$$
\begin{array}{rcll} \qquad
  \u{\a}\u{\a+1}\u{\a}&\equiv&\u{\a+1}\u{\a}\u{\a+1},\\ 
  \u{\b}\u{\c}\u{\a}&\equiv&\u{\b}\u{\a}\u{\c},&\qquad\hbox{if $\a<\b<\c$},\\ 
  \u{\a}\u{\c}\u{\b}&\equiv&\u{\c}\u{\a}\u{\b},&
  \qquad\hbox{if $\a<\b<\c$},\\ 
  \u{\a}\u{\a}&\equiv &0,\end{array}
  \eqno{(2.2)}
$$
the set of all reduced sequences $R(w)$ for a permutation $w\in\Sym_\infty$
is refined into classes, called Coxeter-Knuth cells,
indexed by some semi-standard tableaux.
The cardinality of a cell is the number
of standard tableaux of the same shape as the cell's 
index~\cite{EG,LS82b,Stanley84}.
This decomposition suggests an action of the symmetric group on $R(w)$.
The symmetric function corresponding to such an action is the function $F_w$
introduced by Stanley in~\cite{Stanley84}. 
Equation~(1.3) suggests the possibility of similar
structure for the monoid $\Mon$ and relations (1.1).

\section{$k$-Bruhat orders and the monoid $\Mon$}

The multiplicative structure of Schubert polynomials is determined by 
Monk's rule \cite{Macdonald91}:
$$
{\mathfrak S}_u (x_1+x_2+\cdots+x_k) = 
\sum_{\stackrel{a\le k <b}{\mbox{\scriptsize $\ell(u(a\,b))=\ell(u)+1$}}}
{\mathfrak S}_{u(a\,b)}.
$$
Successive applications of this give
  $$
{\mathfrak S}_u (x_1+x_2+\cdots+x_k)^n = 
\sum_{\stackrel{\mbox{\scriptsize $w\in\Sym_\infty $}}
{\ell(w)=\ell(u)+n}} \gamma(u,w,k) {\mathfrak S}_w, 
$$
where $\gamma(u,w,k)$ counts the sequences of transpositions
$(a_1\,b_1),$ $(a_2\,b_2),$ $\ldots,$ $(a_n\,b_n)$ such that 
$w=u(a_1\,b_1)(a_2\,b_2)\cdots(a_n\,b_n)$ and, for all $r$, we have $a_r\le k<b_r$ with
  $$\ell\big(u (a_1\,b_1)(a_2\,b_2)\cdots(a_{r-1}\,b_{r-1})\big)=
  \ell\big(u (a_1\,b_1)(a_2\,b_2)\cdots(a_r\,b_r)\big)+1.$$
On the other hand 
  $$(x_1+x_2+\cdots+x_k)^n = 
\sum_{\lambda} f^\lambda S_\lambda(x_1,x_2,\ldots,x_k),$$
where $S_\lambda(x_1,x_2,\ldots,x_k)$ is the Schur polynomial indexed by a
partition 
$\lambda$ of $n$. There is a unique Grassmannian permutation
$v(\lambda,k)$ such that 
$S_\lambda(x_1,x_2,\ldots,x_k)=
{\mathfrak S}_{v(\lambda,k)}$~\cite{Macdonald91}.
Hence
 $${\mathfrak S}_u (x_1+x_2+\cdots+x_k)^n = 
   \sum_{\lambda} f^\lambda {\mathfrak S}_u {\mathfrak S}_{v(\lambda,k)} =
   \sum_w\left(\sum_\lambda f^\lambda c_{u\, v(\lambda,k)}^w
   \right){\mathfrak S}_w,
$$
and we have
$$
 \sum_\lambda f^\lambda c_{u\, v(\lambda,k)}^w=\gamma(u,w,k).
 \eqno{(3.1)}
$$

The equation (3.1) suggests that we should study the 
partial order defined by the following
relation: $u\le_k w$ if and only if $\gamma(u,w,k)>0$.
Equivalently, this is the partial order with covering relation given by the
index of summation in Monk's rule.
We call this suborder of the Bruhat order the {\sl $k$-Bruhat order}.
Denote by $[u,w]_k$ the interval from $u$ to $w$ in the $k$-Bruhat order.
Then $\gamma(u,w,k)$ is the number of maximal chains in $[u,w]_k$.

These cover relations give some invariants of the $k$-Bruhat order. For
example, consider the following maximal chain in the 3-Bruhat order:
\silentfootnote{$^{\dag}$Notation: For every $w\in\Sym_\infty$
there exists infinitely many $n$ such that 
$w\in\Sym_n\subset\Sym_\infty$. 
For any such $n$ we write
$(w(1),w(2),\ldots,w(n))$ to represent such a $w$.}
$$
   (3,1,5,2,6,4)\le_3 (3,1,6,2,5,4)\le_3(3,2,6,1,5,4)
    \le_3(3,5,6,1,2,4).^{\dag}
$$
In this chain, the first three entries of the
permutations do not decrease and the other entries do not increase. 
Also, the second and third entries remain in the same relative order for all
permutations in the chain. 
This leads to a characterization of the
$k$-Bruhat order based on such invariants.
\vfill\eject

\begin{prop}[Theorem A of \cite{BS97a}]\label{klength} 
For $u,w\in\Sym_\infty$, \ 
$u\le_k w$ if and only if
\begin{enumerate}
 \item[(1)] $u(i)\le w(i)$ \qquad for $i\le k$,
 \item[(2)] $u(i)\ge w(i)$ \qquad for $i>k$,
 \item[(3)] $(u(i)<u(j)\implies w(i)<w(j))$\quad for\quad 
$i<j\le k$ or $k<i<j$.
\end{enumerate}
\end{prop}

The sufficiency of these conditions follows from the existence of a specific 
maximal chain in the interval $[u,w]_k$. We call it the CM-chain of $[u,w]_k$.

\begin{definition}[CM-chain]\label{canondef} \sl
For $u<_k w$, the CM-chain of the interval $[u,w]_k$ is recursively defined
as follows: 
\begin{enumerate}
\item[$\bullet$] If $\ell(w)=\ell(u)+1$ then the unique chain $u<_kw$ 
is the CM-chain of
$[u,w]_k$.
\item[$\bullet$] If $\ell(w)>\ell(u)+1$, let $a\le k<b$ be the unique
 integers such that 
 \begin{enumerate}
  \item[\bf I] $u(a)<w(a)$ and $w(a)=\max\{w(j):j\le k, u(j)<w(j)\}$,
  \item[\bf II] $u(b)>u(a)\ge w(b)$ and $w(b)=\min\{w(j):j>k, u(j)>u(a)\ge
  w(j)\}$. 
 \end{enumerate}
  Let $u_1=u(a\,b)$. The CM-chain of $[u,w]_k$ is
  $$u=u_0<_ku_1<_ku_2<_k\cdots<_ku_n=w,$$
  where $u_1<_ku_2<_k\cdots<_ku_n$ is the CM-chain of $[u_1,w]_k$.
\end{enumerate}
\end{definition}

It is not obvious that conditions {\bf I} and {\bf II} define unique 
integers $a\le k<b$. 
We refer the reader to \S 3.1 of \cite{BS97a} for a complete proof 
of this fact. 
The symmetry in the conditions (1)-(3) of Proposition~\ref{klength}
implies the following lemma.

\begin{lem}[Vertical Symmetry]\label{vsym} 
Let $m$ be any integer such that $u,w\in\Sym_m$.
Let $\omega_0$ denote the longest element $(m,m-1,\ldots,1)$ of $\Sym_m$. 
Then the map $\Omega_m:\Sym_m\to\Sym_m$ defined by 
$\Omega_m(u)=\omega_0u\omega_0$ is an order preserving
involution. 
That is
 $$u\le_k w \qquad \iff \qquad \Omega_m(u)\le_{m-k} \Omega_m(w).$$
\end{lem}

We use Lemma~\ref{vsym} to define another specific maximal chain in 
the interval $[u,w]_k$.
Given $u,w\in \Sym_m$, apply $\Omega_m$ to the CM-chain of 
$[\Omega_m(u),\Omega_m(w)]_{m-k}$ 
to obtain the DCM-chain of $[u,w]_k$.
We can define it recursively, as in Definition~\ref{canondef}, replacing 
{\bf I} and {\bf II} by:
\begin{itemize}
 \item[\bf I$'$] $u(b)>w(b)$ and $w(b)=\min\{w(j):j> k, u(j)>w(j)\}$,
 \item[\bf II$'$] $u(a)<u(b)\le w(a)$ and 
$w(a)=\max\{w(j):j\le k, u(j)<u(b)\le w(j)\}$.
\end{itemize}

For example, if $u=(2,1,6,4,3,5)$ and $w=(4,5,6,1,2,3)$,
the first step of the procedure for the CM-chain of $[u,w]_3$ 
gives us $(a,b)=(2,4)$.
The full chain is given below, written from bottom to top.
 $$
\begin{array}{ccccc}
(4,5,6,1,2,3)&& (4,5,6,1,2,3)&& (4,5,6,1,2,3)\\ 
(3,5,6,1,2,4)&& (4,3,6,1,2,5)&& (3,5,6,1,2,4)\\ 
(2,5,6,1,3,4)&& (4,1,6,3,2,5)&& (3,4,6,1,2,5)\\ 
(2,4,6,1,3,5)&& (3,1,6,4,2,5)&& (2,4,6,1,3,5)\\ 
(2,1,6,4,3,5)&& (2,1,6,4,3,5)&& (2,1,6,4,3,5)\\ 
           \hbox{CM-Chain}& &\hbox{A Maximal Chain}& &\hbox{DCM-Chain}
\end{array}
$$
Consider a maximal maximal chain of $[u,w]_k$,
$$
  u=u_0<_k u_1 <_k u_2 <_k \cdots <_k u_n=w,
 \eqno{(3.2)}
$$
where $u_{i+1}=u_i(a_i\,b_i)$.
We note that if (3.2) is the CM-chain, then $w(a_i)>w(a_j)$, or
$w(a_i)=w(a_j)$ and $w(b_i)<w(b_j)$ for all $1\le i<j\le n$.
This motivates our definition of inversion. We say that $(i,j)$ is an 
{\sl inversion} of the
chain~(3.2) if $1\le i<j\le n$ and
\begin{itemize}
\item $w(a_i)<w(a_j)$, or
\item $w(a_i)=w(a_j)$ and  $w(b_i)>w(b_j)$.
\end{itemize}

\noindent The inversion set $\mathcal I$ of a chain is the set of all its
inversions. 
In the example above, the inversion set of the middle maximal chain is
$\{(1,2),(1,3),(1,4),(2,3),(2,4)\}$.
$\{1,2\}$.

\begin{lem}\label{canon}
 A maximal chain is the CM-chain if and only if it has no inversions.
\end{lem}

\noindent{\em Proof  }
The reverse implication is clear. Consider a maximal chain with
an inversion $(i,j)$. It suffices to show there is an
$i'$ such that $(i',i'+1)$ is also an inversion of the chain. 
If $(i,i+1)$ is an inversion then we
are done. If $(i,i+1)$ is not an inversion then either
\begin{enumerate}
\item[(a)] $w(a_i)>w(a_{i+1})$ or
\item[(b)] $w(a_i)=w(a_{i+1})$ and $w(b_i)<w(b_{i+1})$.
\end{enumerate}

\noindent
In the first case we have $w(a_{i+1})<w(a_i)\le w(a_j)$, 
and in the second case we have
$w(a_{i+1})=w(a_i)< w(a_j)$, or $w(a_{i+1})=w(a_i)= w(a_j)$ 
and $w(b_{i+1})>w(b_i)> w(b_j)$.
Thus $(i+1,j)$ is an inversion. 
By induction on $j-i$ we conclude that there is an $i\le i'<j$ such 
that $(i',i'+1)$ is an inversion of the chain. 
\qed\medskip

Our next objective is to generate all the maximal chains of $[u,w]_k$.
For this we need the definitions of $\ell_{\u{}}$ and $\preceq$.
The reader will find more details in~\cite{BS97a}.

\begin{prop}[Theorem E of \cite{BS97a}]\label{indep} 
For $u\le_k w$ and $u'\le_{k'} w'$, if
$wu^{-1}=w'(u')^{-1}$, then $v\mapsto vu^{-1}u'$ induces 
$[u,w]_k\cong [u',w']_{k'}$.
\end{prop}

\begin{prop}[Theorem 3.1.5 of \cite{BS97a}]\label{Izeta}  
For $\zeta\in\Sym_\infty$, let
$up(\zeta)=\{j:\zeta^{-1}(j)<j\}=\{j_1<j_2<\cdots<j_k\}$. Let
$w=\big(j_1,j_2,\ldots,j_k,\ldots\big)$ where, to the right of $j_k$, 
we put the complement of
$up(\zeta)$ in increasing order. We have that $[\zeta^{-1}w,w]_k$ is nonempty.
\end{prop}

We use the above two propositions to define the function $\ell_{\u{}}$.
The number $k$ in Proposition~\ref{Izeta} is the smallest possible for which
$[u,w]_k$ is nonempty and $w=\zeta u$. 
The length difference $\ell(w)-\ell(u)$ is
the same for all nonempty $[u,w]_k$ such that $w=\zeta u$. 
With this in mind we define
$\ell_{\u{}}(\zeta)$ to be the length difference 
$\ell(w)-\ell(u)$ obtained from any nonempty
$[u,w]_k$ such that $w=\zeta u$.
This shows part (1) of Theorem~\ref{Operators}. 
Let $dw(\zeta)=\{j:\zeta^{-1}(j)>j\}$.
Proposition~\ref{Izeta} constructs a standard interval $[u,w]_k$ for any
$\zeta$. 
Counting the inversions of $u$ and $w$, and rearranging the terms we deduce
\begin{eqnarray*}
 \ell_{\u{}}(\zeta)&=&
 \abs{\big\{ (i,j)\in up(\zeta)\times dw(\zeta) : \,i>j\big\}}\\
 &&-\abs{\big\{ (i,j)\in \zeta^{-1}(up(\zeta))\times
 \zeta^{-1}(dw(\zeta)):\,i>j\big\}}\\
 &&-\abs{\big\{ (i,j)\in \left(\zeta^{-1}(up(\zeta))\right)^2:\,i<j 
 \hbox{ and }  \zeta(i)>\zeta(j) \big\}}\\
 &&-\abs{\big\{ (i,j)\in \left(\zeta^{-1}(dw(\zeta))\right)^2:\,i<j 
 \hbox{ and }\zeta(i)>\zeta(j) \big\}}.
\end{eqnarray*}

If $u\le_k \zeta u$ and $u'\le_{k'} \zeta u'$, then the isomorphism
$[u,\zeta u]_k\cong [u',\zeta u']_{k'}$ is given by 
$\eta u\mapsto \eta u'$. We now
introduce the universal $k$-Bruhat order $\preceq$ on $\Sym_\infty$.
$$
  \eta\preceq\zeta \quad\iff\quad\big(\hbox{There exists $u$, $k$ such that
  $u\le_k\eta u\le_k\zeta u$}\big).
  \eqno{(3.3)}
$$
Using the permutation $u$ given by Proposition ~\ref{Izeta}, 
we see that $\eta\preceq \zeta$ if
and only if
\begin{enumerate}
 \item[(1)] $\a\le\eta(\a)\le \zeta(\a)$ \quad for 
  $\a\in\zeta^{-1}\big(up(\zeta)\big)$,
 \item[(2)] $\a\ge\eta(\a)\ge \zeta(\a)$ \quad for 
  $\a\in\zeta^{-1}\big(dw(\zeta)\big)$,
 \item[(3)] $\big(\eta(\a)<\eta(\b)\implies \zeta(\a)<\zeta(\b)\big)$ for
  $\a<\b\in \zeta^{-1}(up(\zeta))$ or $\a<\b\in \zeta^{-1}(dw(\zeta))$.
\end{enumerate}

\noindent It follows from the definition that the order $\preceq$ is ranked by
$\ell_{\u{}}$ and $[1,\zeta\eta^{-1}]_\preceq\cong[\eta,\zeta]_\preceq$ 
via the map
$\xi\mapsto\xi\eta$. 
The operators $\hu{\a\b}$ in~(1.2) are defined so that
$\hu{\a\b}\eta=\zeta$ if and only if $\zeta$ covers $\eta$ in $\preceq$. 
In particular, nonzero compositions 
$\hx=\hu{\a_n\b_n}\cdots\hu{\a_2\b_2}\hu{\a_1\b_1}$ such that 
$\hx\eta=\zeta$ correspond bijectively to maximal chains in 
$[\eta,\zeta]_\preceq$:
$$
 \eta\preceq \hu{\a_1\b_1}\eta\preceq \hu{\a_2\b_2}\hu{\a_1\b_1}
 \eta\preceq\cdots\preceq\hx\eta=\zeta
$$
We note that the isomorphism 
$[1,\zeta\eta^{-1}]_\preceq\cong[\eta,\zeta]_\preceq$ implies
$$
 \hx \eta=\zeta\quad\iff\quad \hx 1=\zeta\eta^{-1}.
 \eqno{(3.4)}
$$

The isomorphism $[1,wu^{-1}]_\preceq\cong[u,w]_k$ given by 
$\eta\mapsto \eta u$,
induces an isomorphism on chains. 
Given a maximal chain
$$
 u=u_0<_k u_1 <_k u_2 <_k \cdots <_k u_n=w 
 \eqno{(3.5)}
$$
of $[u,w]_k$, we adopt the following conventions.
\begin{itemize}
 \item Let $a_i\le k<b_i$ be such that $u_{i+1}=u_i(a_i\,b_i)$.
 \item Let $\a_i=u_{i-1}(a_i)$ and $\b_i=u_{i-1}(b_i)$. Hence
  $u_i=(\a_i\,\b_i)u_{i-1}$.
\end{itemize}

\noindent Under the isomorphism above, this defines a unique (nonzero)
composition
$$
 \hx=\hu{\a_n\b_n}\cdots\hu{\a_2\b_2}\hu{\a_1\b_1}
 \eqno{(3.6)}
$$  
such that $wu^{-1}=\hx 1$. 
Conversely, given a nonzero composition as in~(3.6) such that 
$wu^{-1}=\hx 1$, we define a
unique maximal chain as in~(3.5) where 
$u_i=(\hu{\a_i\b_i}\cdots\hu{\a_1\b_1} 1)u$. 
This correspondence is used to encode maximal chains for the rest of the 
paper. 
Via this identification, we will refer to a nonzero composition 
$\hx$ such that $\hx 1=wu^{-1}$ as a maximal chain of $[u,w]_k$.

Proposition~\ref{Izeta} is very useful for constructing 
intervals in $k$-Bruhat orders.
For example, let
$\zeta=(5,4,2,1,3)$. Proposition~\ref{Izeta} gives 
$u=(2,1,4,3,5)\le_2(4,5,1,2,3)=\zeta u$. 
{}From Definition~\ref{canondef}, the CM-chain is
$\hu{34}\hu{23}\hu{45}\hu{14}$. 
Now if we apply the relations (1)-(3) of (1.1) to the
CM-chain  we get:
$$
\hu{34}\underline{\hu{23}\hu{45}}\hu{14}
  \equiv\hu{34}\hu{45}\underline{\hu{23}\hu{14}}
  \equiv\underline{\hu{34}\hu{45}\hu{14}}\hu{23}
  \equiv\hu{35}\underline{\hu{13}\hu{34}\hu{23}}
  \equiv\hu{35}\hu{23}\hu{12}\hu{24}.
$$ 
These are all the maximal chains in the interval $[u,w]_k$ as depicted in 
Figure 1.
The first two equivalences are instances of the relation (3) 
of~(1.1), the last two are
instances of relations (1) and (2) of (1.1), respectively. 
The second chain is the DCM-chain. 

\begin{figure}[htb]
$$
\epsfxsize=3.in \epsfbox{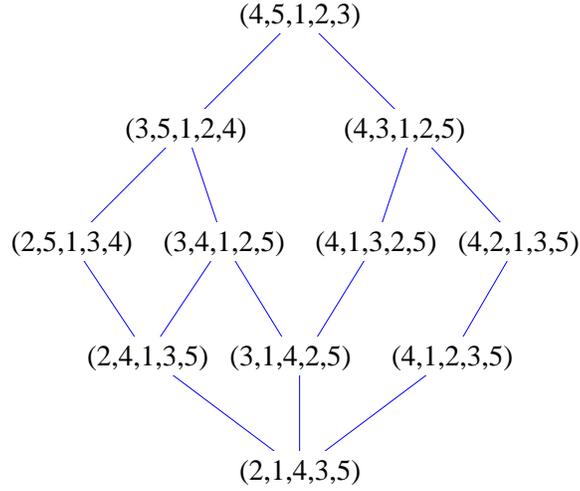}
$$
\caption{The interval $[(2,1,4,3,5),\: (4,5,1,2,3)]_2$.
\label{fig:one}}
\end{figure}

\begin{thm}\label{genere} 
If $u\le_k w$, then any two maximal chains in $[u,w]_k$ are connected
by a series of relations (1)-(3) of (1.1). 
Moreover, it is never possible to apply any of
the relations (4) or (5) of (1.1) to a maximal chain.
\end{thm}

\noindent{\em Proof  }
We first show that any of the relations (1)-(3) of~(1.1) 
that can be applied to a maximal chain
$$   \hx =\hu{\a_n\b_n}\cdots\hu{\a_2\b_2}\hu{\a_1\b_1} 
   \eqno{(3.7)}
$$
in $[u,w]_k$ results in another maximal chain. 
Moreover, the relations (4) and (5) can never
be applied to this chain. Given the maximal chain~(3.7), let
$u_i=(\hu{\a_i\b_i}\cdots\hu{\a_1\b_1} 1)u$ be as before, for $0\le i\le n$. 
Then since
$u_{i-1}\le_k u_i$ is a cover,
\begin{enumerate}
\item[(i)] $u_i=(\alpha_i\,\,\beta_i)u_{i-1}=u_{i-1}(a_i\,\,b_i)$ with 
$a_i\le k<b_i$.
\item[(ii)] If $\alpha_i<\gamma<\beta_i$, then  
$u^{-1}_{i-1}(\gamma)<a_i$ or $b_i<u^{-1}_{i-1}(\gamma)$. 
\end{enumerate}

Consider applying the relations~(1.1) to a segment of length 
two in the chain~(3.7).
We may assume that the segment is $\hu{\a_2\b_2}\hu{\a_1\b_1}$.
Suppose $\{\a_1,\b_1\}\cap\{\a_2,\b_2\}=\emptyset$, and assume $\a_1<\a_2$, 
as the other case
is symmetric. There are three possible relative orders for the numbers
$\a_1,\b_1,\a_2$ and $\b_2$. We consider each in turn. 
If $\a_1<\a_2<\b_1<\b_2$,
the situation in relation (4) with strict inequalities,
then condition (ii) for $i=1$ implies $a_2=u_0^{-1}(\a_2)<a_1$, and for 
$i=2$ implies
$a_1=u_1^{-1}(\b_1)<a_2$, a contradiction. 
Now suppose $\a_1<\b_1<\a_2<\b_2$ or $\a_1<\a_2<\b_2<\b_1$. 
An example of each case is found as a square in Figure 1. 
Then (i) and (ii) impose no additional conditions on $a_1,a_2,b_1$ and
$b_2$, so 
$u_0\le_k u_0(a_2\,\, b_2)\le_k u_0(a_2\,\, b_2)(a_1\,\,b_1)=u_2$.

Suppose one of the relations (1) or (2) of~(1.1) applies to a 
segment of length three.
Again an example of each case is found as a hexagon in Figure 1.
Both argument are similar, so suppose that (1) applies.
We have $\a<\b<\c<\d$ and the segment is $\hu{\b\c}\hu{\c\d}\hu{\a\c}$.
By condition (ii), the numbers $\a,\b,\c$ and $\d$ appear in $u$ 
in one of the following two
orders
$$
\big(\ldots,\b,\ldots,\a,\ldots,\c,\ldots,\d,\ldots\big) 
\quad \hbox{ or } \quad
\big(\ldots,\b,\ldots,\a,\ldots,\d,\ldots,\c,\ldots\big).
$$
Suppose we are in the first case, the argument in the second being similar. 
Then the chain is
  $$
\begin{array}{c}
  \big(\ldots,\c,\ldots,\d,\ldots,\a,\ldots,\b,\ldots\big)\\ 
  \big(\ldots,\b,\ldots,\d,\ldots,\a,\ldots,\c,\ldots\big)\\ 
  \big(\ldots,\b,\ldots,\c,\ldots,\a,\ldots,\d,\ldots\big)\\ 
  \big(\ldots,\b,\ldots,\a,\ldots,\c,\ldots,\d,\ldots\big)
\end{array}
$$
It is clear that
  $$
\begin{array}{c}
  \big(\ldots,\c,\ldots,\d,\ldots,\a,\ldots,\b,\ldots\big)\\ 
  \big(\ldots,\c,\ldots,\b,\ldots,\a,\ldots,\d,\ldots\big)\\ 
  \big(\ldots,\c,\ldots,\a,\ldots,\b,\ldots,\d,\ldots\big)\\ 
  \big(\ldots,\b,\ldots,\a,\ldots,\c,\ldots,\d,\ldots\big)\\ 
\end{array}
$$
is also a chain. This is represented by $\hu{\b\d}\hu{\a\b}\hu{\b\c}$, 
completing this case. 
To conclude our first objective,
we notice that the fourth relation, with equalities, or the fifth relation, 
are clearly not possible for $k$-Bruhat orders, by 
Proposition~\ref{klength} (1) and (2).

We now show that any two maximal chains in $[u,w]_k$ are connected by 
successive uses of the relations~(1.1). 
It suffices to show that any maximal chain $\hx$ is connected to
the CM-chain. 
For this we proceed by induction on $n$. 
If $n=1$, then there is a unique maximal chain. 
Let $n>1$ and assume that the theorem holds for all intervals 
$[u',w']_{k'}$ such that
$\ell(w')-\ell(u')<n$. 
That is, we may assume that $\hx= \hy\hu{\a_1\b_1}$ where $\bf y$ is
any maximal chain. 
If $a_1,b_1$ satisfy the conditions 
{\bf I} and {\bf II} of Definition~\ref{canondef} then
choosing $\hy$ to be the CM-chain of $[u_1,w]_k$ completes the proof 
since then $\hx$ is the
CM-chain of $[u,w]_k$. 
If condition {\bf I} fails, then $w(a_1)$ is not maximal with
$u(a_1)<w(a_1)$. 
In this case assume that $\hy$ is the CM-chain of $[u_1,w]$ so that
$w(a_2)>w(a_1)$. 
We have two sub-cases to consider:

\noindent Case 1a: 
$\{\a_1,\b_1\}\cap\{\a_2,\b_2\}=\emptyset$. We can use relation (3) of
(1.1) and get
$$
 \hx\equiv \hu{\a_n\b_n}\cdots\hu{\a_1\b_1}\hu{\a_2\b_2}.
 \eqno{(3.8)}
$$
The hypothesis on $\bf y$ and $w(a_2)>w(a_1)$ implies that 
$\hu{\a_2\b_2}$ is the first step of
the CM-chain of $[u,w]_k$. 
We can use our induction hypothesis on
$[\hu{\a_2\b_2}u,w]_k$ and get $\hx\equiv \hz\hu{\a_2\b_2}$, t
he CM-chain of $[u,w]_k$.

\noindent Case 1b: 
$\a_2<\b_2=\a_1<\b_1$. Since $\bf y$ is the CM-chain of $[u_1,w]_k$, we have
 $$
\b_2=\a_3<\b_3=\a_4<\cdots<\b_{m-1}=\a_m,
$$ 
for $m\ge 3$, where $\b_m =
w(a_2)>w(a_1)\ge\b_1$. 
Let $3\le s\le m$ be such that $\a_s<\b_1<\b_s$. We can apply the relations
(1.1) and get
\begin{eqnarray*}
\hx&=&\hu{\a_n\b_n}\cdots \hu{\a_m\b_m}\cdots\hu{\a_s\a_s}\
         cdots\hu{\a_2\b_2}\hu{\a_1\b_1}\\ 
     &\equiv&\hu{\a_n\b_n}\cdots \hu{\a_m\b_m}\cdots\hu{\a_{s+1}\a_{s+1}}
             \hu{\a_s\b_1}\hu{\b_1\b_s}\hu{\a_2\b_1}\hu{\a_{s-1}\b_{s-1}}
              \cdots\hu{\a_3\b_3}\\ 
     &\equiv&\hu{\a_n\b_n}\cdots \hu{\a_m\b_m}\cdots\hu{\a_{s+1}\a_{s+1}}
             \hu{\a_s\b_1}\hu{\a_{s-1}\b_{s-1}}\cdots\hu{\a_3\b_3}
              \hu{\b_1\b_s}\hu{\a_2\b_1}\\ 
     &\equiv&\hz\hu{\a_2\b_1}.
\end{eqnarray*}
where, by the induction hypothesis, $\hz$ is the CM-chain of
$[\hu{\a_1\b_2}u,w]_k$. 
Here $\hu{\a_2\b_1}$ is the first step in the CM-chain of
$[u,w]_k$. Hence
$\hx\equiv \hz\hu{\a_2\b_1}$, the CM-chain of $[u,w]_k$.

If condition {\bf I} holds but condition {\bf II} fails, 
then $w(b_1)$ is not minimal.
In this case assume that $\bf y$ is the DCM-chain of $[u_1,w]$. 
Here, we must have that
$w(b_2)<w(b_1)$ and again we have two sub-cases to consider:

\noindent Case 2a: 
$\{\a_1,\b_1\}\cap\{\a_2,\b_2\}=\emptyset$. We can use the relation (3) of
(1.1) and the induction hypothesis to get
$$
 \hx\equiv \hu{\a_n\b_n}\cdots\hu{\a_1\b_1}\hu{\a_2\b_2}
 \equiv\hz\hu{\a_2\b_2},
 \eqno{(3.9)}
$$
where $\hz$ is the CM-chain of $[\hu{\a_2\b_2}u,w]_k$. 
If $\hu{\a_2\b_2}$ is the first step in
the CM-chain of $[u,w]_k$ we are done. 
If not, then condition {\bf I$'$} on $\hu{\a_2\b_2}$
implies that only condition {\bf I} can fail in
$\hz\hu{\a_2\b_2}$ and we are back to cases 1a or 1b. 

\noindent Case 2b: 
$\a_1<\b_1=\a_2<\b_2$. Since $\bf y$ is the DCM-chain of $[u_1,w]_k$, we have
 $$\a_2=\b_3>\a_3=\b_4>\cdots>\a_{m-1}=\b_m,$$ 
for $m\ge 3$, where $\a_m =
w(b_2)>w(b_1)\ge\a_1$. Let $3\le s\le m$ be such that $\b_s>\a_1>\a_s$. 
We can apply the relations~(1.1) and get
$$
 \begin{array}{rcl}
   \hx&=&\hu{\a_n\b_n}\cdots \hu{\a_m\b_m}\cdots\hu{\a_s\a_s}\cdots
          \hu{\a_2\b_2}\hu{\a_1\b_1}\\ 
      &\equiv&\hu{\a_n\b_n}\cdots \hu{\a_m\b_m}\cdots\hu{\a_{s+1}\a_{s+1}}
             \hu{\a_1\b_s}\hu{\a_s\b_1}\hu{\a_1\b_2}\hu{\a_{s-1}\b_{s-1}}
             \cdots\hu{\a_3\b_3}\\ 
      &\equiv&\hu{\a_n\b_n}\cdots \hu{\a_m\b_m}\cdots\hu{\a_{s+1}\a_{s+1}}
             \hu{\a_1\b_s}\hu{\a_{s-1}\b_{s-1}}\cdots\hu{\a_3\b_3}
             \hu{\a_s\b_1}\hu{\a_1\b_2}\\ 
      &\equiv& \hz\hu{\a_1\b_2},
 \end{array}
 \eqno{(3.10)}
$$
where $\hz$ is the CM-chain of $[\hu{\a_1\b_2}u,w]_k$. 
If $\hu{\a_1\b_2}$ is the first step in the
CM-chain of $[u,w]_k$, then we are done. If not, then condition {\bf I$'$} on
$\hu{\a_1\b_2}$ implies that only the condition {\bf I} can fail in
$\hz\hu{\a_2\b_2}$ and again we are back to cases 1a or 1b. 
\qed\medskip

We now complete the characterization of compositions 
$\hx=\hu{\a_n\b_n}\cdots\hu{\a_1\b_1}$
which correspond to maximal chains for some $[u,w]_k$. 
If $\hx$
corresponds to a maximal chain in $[u,w]_k$, then
$wu^{-1}=\hx$. Hence $w=\zeta u$ where
$\zeta=\hx 1=wu^{-1}$. Conversely, Proposition~\ref{Izeta} shows that
for any $\zeta\in\Sym_\infty$ we can find $u$ and $w$ such that 
$w=\zeta u$ and $[u,w]_k$ is
nonempty for some $k$.
In the following, we say that a composition
$\hx=\hu{\a_n\b_n}\cdots\hu{\a_1\b_1}$ is $\u{}$-reduced if $\hx 1\neq 0$.
Theorem~\ref{genere} gives us a way of generating all
$\u{}$-reduced sequences for $\zeta\in\Sym_\infty$; 
they are all connected via the
relations (1)-(3) of (1.1). 
To complete our study, we need to characterize the
compositions $\hx$ such that $\hx =0$.

\begin{thm}\label{reduced} 
Let $\hx=\hu{\a_n\b_n}\cdots\hu{\a_1\b_1}$ be a composition.
If $\hx 1=0$, then $\hx\equiv 0$ modulo the relations (1.1).
\end{thm}

\noindent{\em Proof  }
We proceed by induction on $n$. 
When $n=2$, \  $\hx 1=0$ implies that relation (4) applies to
$\hx$. 
Suppose $n\ge 3$ and the theorem holds for all compositions of length $<n$.
Let $\hy=\hu{\a_{n-1}\b_{n-1}}\cdots\hu{\a_1\b_1}$ and 
we may assume that $\hy 1=\tau\neq 0$.

We first characterize those $w$ such that $\tau^{-1}w\le_k w$, for some $k$.
Let $up(\tau)$ and $dw(\tau)$ be defined as above, and let $fix(\tau)$ be 
the set of fixed points of $\tau$. By Proposition~\ref{klength}, 
$u=\tau^{-1}w\le_k w$ if and only if
$$
 \begin{array}{l}
      \bullet\quad up(\tau)\subseteq 
 \{w(i):1\le i\le k\}\subseteq up(\tau)\cup fix(\tau),\\ 
 \bullet\quad \hbox{for $i<j\le k$ or $k<i<j$, 
 if $u(i)<u(j)$ then $w(i)<w(j)$.}
 \end{array}
 \eqno{(3.11)}
$$
The second condition implies that if $\a<\c$ are in 
$up(\tau)\cup fix(\tau)$ and
$\tau^{-1}(\a)>\tau^{-1}(\c)$, then 
$\max\{w^{-1}(\a),w^{-1}(\c)\}\le k$ implies $w^{-1}(\a)<w^{-1}(\c)$.
Similarly, if $\c<\b$ are in $dw(\tau)\cup fix(\tau)$ and
$\tau^{-1}(\c)>\tau^{-1}(\b)$ then $k<\min\{w^{-1}(\b),w^{-1}(\c)\le k\}$ 
implies $w^{-1}(\c)<w^{-1}(\b)$.
With this and the definition of $\preceq$, we see that
$\ell_{\u{}}\big((\a_n,\b_n)\tau\big)\neq\ell_{\u{}}(\tau)+1$ 
implies one of the following
holds:
\begin{enumerate}
\item[(a)] $\a_n\in dw(\tau)$,
\item[(b)] $\b_n\in up(\tau)$,
\item[(c)] $\a_n<\gamma<\b_n$ where $\tau^{-1}(\a_n)>\tau^{-1}(\gamma)$, or
$\tau^{-1}(\gamma)>\tau^{-1}(\b_n)$.
\end{enumerate}

\noindent
We complete the proof by showing that each case (a), (b), or (c) implies
$\hu{\a_n\b_n}\hy\equiv 0$ modulo the relations (1.1).

\noindent{\sl If (a) holds:} 
By Theorem~\ref{genere} we may assume that $\hy$ is any maximal chain.
Let $\hy=\hz\hu{\a_1\b_1}$.
Note that if $\a_n\in dw\big(\hz 1)$ then
the induction hypothesis applies and we are done. We can thus assume that
$\a_n=\a_1$. But this must be true for any maximal chain $\hy$.
Since $\a_1=\min\big(dw(\tau)\big)$
for the DCM-chain, we have $\a_n=\min\big(dw(\tau)\big)$.
Now let $\hy$ be the CM-chain, and consider its initial segment 
$\hu{\a_m\b_m}\cdots\hu{\a_1\b_1}$ where 
$\b_1=\a_2<\b_2=\a_3<\cdots<\b_{m-1}=\a_m$ 
and $\b_m =\max\big(up(\tau)\big)$.
If $\big|up(\tau)\big|>1$, then  $m<n-1$. 
Consider the next operator $\hu{\a_{m+1}\b_{m+1}}$.
Since $\a_1=\min\big(dw(\tau)\big)$, we have $\a_1<\a_m+1$, and since 
$\b_m=\max\big(up(\tau)\big)$, we have $\b_{m+1}<\b_m$.
Thus we may apply a sequence of the relations (1)-(3) of (1.1), 
as in~(3.10), to
obtain $\hy\equiv \hz'\hu{\a_{m+1}\b'}$
for some $\hz'$ and $\b'$. Since $\a_n=\a_1\in dw\big(\hz' 1)$, the induction
hypothesis applies to conclude $\hu{\a_n\b_n}\hz'\equiv 0$. 
Thus we may assume that (a) holds and
$\big|up(\tau)\big|=1$. That is,
$\b_1=\a_2<\b_2=\a_3<\cdots<\b_{n-2}=\a_{n-1}$ and $\a_n=\a_1$. 
If $\b_n<\a_{n-1}$ or
$\b_n>\b_{n-1}$ then we apply relation (3) to obtain
$\hu{\a_n\b_n}\hy\equiv
\hu{\a_{n-1}\b_{n-1}}\hu{\a_n\b_n}\hu{\a_{n-2}\b_{n-2}}\cdots\hu{\a_1\b_1}
    \equiv \hu{\a_{n-1}\b_{n-1}}\hy'$, 
and $\hy'\equiv 0$ by the induction hypothesis.
If $\b_n=\a_{n-1}$ then we may apply relation (2) to obtain 
$\hu{\a_n\b_n}\hx\equiv
    \hu{\a_{n-2}\b_{n-2}}\hu{\a_n\a_{n-2}}
\hu{\a_{n-2}\b_{n-1}}\cdots\hu{\a_1\b_1}$, which is
equivalent to $0$ as before. Finally if $\a_{n-1}<\b_n\le\b_{n-1}$ then
$\hu{\a_n\b_n}\hu{\a_{n-1}\b_{n-1}}\equiv 0$

\noindent{\sl If (b) holds:} 
This case is similar to (a), the map $\Omega_n$ from Lemma~\ref{vsym}
can be used to interchange the roles of conditions (a) and (b).

\noindent{\sl If (c) holds:}
Assume that $\tau^{-1}(\a_n)>\tau^{-1}(\c)$.
The other case, $\tau^{-1}(\c)>\tau^{-1}(\b_n)$, 
is argued in a similar fashion using the map
$\Omega_n$. 
We may also assume that (a) does not hold, hence we have
$\tau^{-1}(\c)<\tau^{-1}(\a_n)\le
\a_n<\c<\b_n$ and, in particular, $\c\in up(\tau)$.
Let $\c$ be minimal with these properties.
We may assume that $\hy$ is the CM-chain and we let 
$\hy=\hu{\a_{n-1}\b_{n-1}}\hz$ and
$\hz 1=\sigma\in\Sym_\infty$. In this case
$\b_{n-1}=\min\big(up(\tau)\big)\le\c$. 
If $\b_{n-1}<\c$ then the minimality of $\c$ implies
$\b_{n-1}\le\a_n$.  
We have a four sub-cases to consider:
\begin{enumerate}
\item[(i)] If $\b_{n-1}=\c$ and $\a_{n-1}\le\a_n$, then 
$\hu{\a_n\b_n}\hu{\a_{n-1}\b_{n-1}}\equiv 0$ is an instance of 
relation (4) of (1.1).
\item[(ii)] If $\b_{n-1}=\c$ and $\a_{n-1}>\a_n$, then 
 $\hu{\a_n\b_n}\hy
   \equiv \hu{\a_{n-1}\b_{n-1}}\hu{\a_n\b_n}\hz.$
Since $\tau^{-1}(\c)<\a_n$ and $\bf x$ is the CM-chain, we must
have $\b_{n-2}=\a_{n-1}$. So
$\a_n<\a_{n-1}=\b_{n-2}<\c<\b_n$ and 
$\sigma^{-1}(\b_{n-2})=\tau^{-1}(\c)<\a_n$. 
By the induction hypothesis $\hu{\a_n\b_n}\hy\equiv 0$.
\item[(iii)]
If $\b_{n-1}<\a_n$, then $\hu{\a_n\b_n}\hy
   \equiv \hu{\a_{n-1}\b_{n-1}}\hu{\a_n\b_n}\hz$ where 
$\sigma^{-1}(\b_{n-2})=\tau^{-1}(\c)$.
The induction hypothesis applies and again  $\hu{\a_n\b_n}\hz\equiv 0$.
\item[(iv)]
If $\b_{n-1}=\a_n$, then since $\hy$ is the CM-chain, the 
minimality of $\c$ implies that
$\b_m=\c<\b_n$ for some $1\le m\le n-2$, with
 $$\a_n=\b_{n-1}>\a_{n-1}=\b_{n-2}>\cdots>\a_{m+2}=\b_{m+1}>\a_{m+1}.$$
For some $1\le s\le m$ we also have
 $$\c=\b_{m}>\a_m=\b_{m-1}>\cdots>\a_{s+1}=\b_s,$$
where $\a_s=\tau^{-1}(\c)<\tau^{-1}(\a_n)=\a_{m+1}$. If
$s>1$ we may appeal to the induction hypothesis and get 
$\hu{\a_n\b_n}\hy\equiv 0$.
Thus we may assume that $s=1$. 
Also, since $\a_1<\a_{m+1}<\b_{m+1}\le\a_n<\c=\b_m$ we may
apply relations (1)-(3) as in~(3.10) to obtain
 $$\hu{\a_{m+1}\b_{m+1}}\hu{\a_{m}\b_{m}}\cdots\hu{\a_1\b_1}\equiv
    \hu{\a'_{m}\b'_{m}}\cdots\hu{\a'_1\b'_1}\hu{\a_{m+1}\b_{m+1}},$$
where $\c=\b_m=\b'_m$, $\b'_{m-1}=\a'_m$, $\b'_{m-2}=\a'_{m-1}$, 
$\ldots$, $\b'_1=\a'_2$ and
$\a'_1=\a_1$. Hence we can use the induction hypothesis on 
$\hu{\a_n\b_n}\cdots\hu{\a_{m+2}\b_{m+2}}
\hu{\a'_{m}\b'_{m}}\cdots\hu{\a'_1\b'_1}$, to obtain
$$\hu{\a_n\b_n}\cdots\hu{\a_{m+2}\b_{m+2}}\hu{\a'_{m}\b'_{m}}
\cdots\hu{\a'_1\b'_1}\equiv 0,$$
and this concludes our proof.
\end{enumerate}
\qed\medskip

\noindent{\em Proof  } [of Theorem \ref{Operators}]
\begin{enumerate}
\item[(a)] This is a direct consequence of Proposition~\ref{indep} 
and Proposition~\ref{Izeta}.
\item[(b)] Theorem~\ref{genere} and Theorem~\ref{reduced} imply 
that the operators $\hu{\a\b}$ satisfy
the relations (1.1). Equation~(3.4) gives the 
characterization part.
\item[(c)] This is a consequence of (b), Theorem~\ref{genere}, 
and Theorem~\ref{reduced}.
\item[(d)] Injection is from part (b) and (c). 
Surjection is given by Proposition~\ref{Izeta}.
\item[(e)] Follows from the definitions of $\preceq$ and $\hu{\a\b}$.
\item[(f)] This is a direct
consequence (a)-(f) above.
\end{enumerate}
\qed\medskip

The universal $k$-Bruhat order is a very interesting object to study on its
own. 
Numerous other results of \cite{BS97a}  can be translated to the monoid
$\Mon$ and on the universal $k$-Bruhat order. We consider a few in the next
section.

\section{A Combinatorial description of $c_\lambda^\zeta$.}

We give a combinatorial description of the constants $c_\lambda^\zeta$
appearing in Equation~(1.3) and combinatorial proofs of many of the
identities of~\cite{BS97a}.  
Recall that the Schur polynomial $S_\lambda(x_1,x_2,\ldots,x_k)$ equals
${\mathfrak S}_{v(\lambda,k)}$ for a unique Grassmannian 
permutation $v(\lambda,k)$.
We have
$$
 {\mathfrak S}_{u}{\mathfrak S}_{v(\lambda,k)}=\sum_{w} 
 c_{u v(\lambda,k)}^w {\mathfrak S}_{w}.
 \eqno{(4.1)}
$$
First we consider a special case of~(4.1).
The Schubert polynomial
${\mathfrak S}_{v((n),k)}=$\break $h_n(x_1,x_2,\ldots,x_k)$ is the
homogeneous symmetric polynomial 
on
$k$ variables. 
Lascoux and Sch\"utzenberger~\cite{LS82a} formulated a Pieri-type formula
for ${\mathfrak S}_u{\mathfrak S}_{v((n),k)}$. 
In~\cite{BB}, proven in~\cite{Sottile96}, we have reformulated this rule.
Using Theorem \ref{Operators}, we can state it here as follows:  
$$
 {\mathfrak S}_u{\mathfrak S}_{v((n),k)} = \sum_
      {\stackrel{\mbox{\scriptsize 
 $\hx=\hu{\a_n\b_n}\cdots \hu{\a_1\b_1}\not\equiv 0 $}}
 {\a_1<\a_2<\cdots <\a_n}}
        {\mathfrak S}_{(\hx 1)u}. 
 \eqno{(4.2)}
$$
There are now other proofs of~(4.2), some of them are 
combinatorial~\cite{Postnikov,Veigneau}.
Let $p=(p_1,p_2,\ldots,p_r)$ be a sequence of $r$ integers such 
that $p_1+p_2+\cdots+p_r=n$.
We say  that a $\u{}$-composition 
$\hx=\hu{\a_n\b_n}\cdots\hu{\a_1\b_1}$ {\sl weakly fits} 
$p$ if 
 $$
\begin{array}{r}
    \hfill\a_1<\a_2<\cdots<\a_{p_1},\\ 
    \hfill\a_{p_1+1}<\a_{p_1+2}<\cdots<\a_{p_1+p_2},\\  
    \hfill\vdots\;\;\\  
    \hfill\a_{n-p_{r}+1}<\a_{n-p_{r}+2}<\cdots<\a_n,
\end{array}
$$ 
and for all $i$, we have $p_i\ge 0$.
Let ${\mathcal H}_p(\zeta)=
\{\hx\in R_{\u{}}(\zeta) : 
\zeta=\hx1\ \hbox{\sl and $\hx$ weakly fits }\,p\}$. 
Note that ${\mathcal H}_p(\zeta)=\emptyset$ if some $p_i<0$.

\begin{Remark}\label{HPier} 
{}From~(4.2), ${\mathcal H}_p(wu^{-1})$ is the 
coefficient of ${\mathfrak S}_w$ in the product 
$$
  {\mathfrak S}_u{\mathfrak S}_{v((p_1),k)}{\mathfrak S}_{v((p_2),k)}\cdots
  {\mathfrak S}_{v((p_r),k)}
$$ 
when all $p_i>0$.
\end{Remark}

\noindent
Now consider the Jacobi identity~\cite{Macdonald91}: for
$\lambda=(\lambda_1,\lambda_2,\ldots,\lambda_r)$ a partition of $n$,
$$
{\mathfrak S}_{v(\lambda,k)}=S_\lambda(x_1,x_2,\ldots,x_k)=
      \det\left(h_{\lambda_i+j-i}(x_1,x_2,\ldots,x_k)
      \right)_{1\le i,j\le r},
 \eqno{(4.3)}
$$
where $h_0(x_1,x_2,\ldots,x_k)=1$, 
$h_n(x_1,x_2,\ldots,x_k)={\mathfrak S}_{v((n),k)}$ for $n>0$, and
$h_n=0$ for $n<0$.  
For $\sigma\in\Sym_r$,
let  $\lambda_\sigma=\big(\lambda_\sigma(1),\lambda_\sigma(2),
\ldots,\lambda_\sigma(r)\big)$,
where
$\lambda_\sigma(i)=\lambda_{\sigma(i)}+i-\sigma(i)$. 
Denote by $\epsilon(\sigma)$ the sign of the
permutation
$\sigma\in\Sym_r$. Expanding the determinant~(4.3) 
in~(4.1), and using~(4.2), we get
\begin{eqnarray*}
   {\mathfrak S}_u{\mathfrak S}_{v(\lambda,k)} 
      &=& \sum_{\sigma\in\Sym_r} \epsilon(\sigma)
   {\mathfrak S}_u{\mathfrak S}_{v((\lambda_\sigma(1)),k)}
      {\mathfrak S}_{v((\lambda_\sigma(2)),k)}\cdots
      {\mathfrak S}_{v((\lambda_\sigma(r)),k)}\\ 
   &=& \sum_{w\in\Sym_\infty}\left( \sum_{\sigma\in\Sym_r} 
       \epsilon(\sigma)\left|
      {\mathcal H}_{\lambda_\sigma}(wu^{-1})\right|\right) {\mathfrak S}_w.
\end{eqnarray*}
Thus
$$
 c_{uv(\lambda,k)}^w=\sum_{\sigma\in\Sym_r} \epsilon(\sigma) \left|
   {\mathcal H}_{\lambda_\sigma}(wu^{-1})\right|.
 \eqno{(4.4)}
$$
This is a consequence of Theorem \ref{Operators}.  
{}From this we deduce the following proposition.

\begin{prop}\label{XX}\mbox{ }
\begin{enumerate}
\item[(1)] $c_{uv(\lambda,k)}^w=0$ if $u\not\le_k w$, and
\item[(2)] if $u\le_k w$ then $c_{uv(\lambda,k)}^w$ depends only on 
$\lambda$ and $wu^{-1}$.
\end{enumerate}
\end{prop}

\noindent Hence, we have that $c_\lambda^\zeta=c_{uv(\lambda,k)}^w$ is 
well defined for any
$u\le_k w$ with $\zeta=wu^{-1}$. We have

\begin{thm}\label{combiconst} 
$\displaystyle c_\lambda^\zeta=\sum_{\sigma\in\Sym_r}
\epsilon(\sigma)
\left| {\mathcal H}_{\lambda_\sigma}(\zeta)\right|$.
\end{thm}

Let us illustrate Theorem~\ref{combiconst} on an example. 
Let $\zeta=(2,5,4,1,6,3)$. 
Using Proposition~\ref{Izeta} we have 
$(3,1,2,5,6,4)=u \le_4 \zeta u = (4,2,5,6,3,1)$. 
In Figure 2, we have drawn the
interval $[u,\zeta u]_4$ and we
have labeled each covering edge in the interval by the index 
$\a$ of the corresponding $\hu{\a\b}$. 
Here we have removed the commas and parentheses to represent the 
permutations in a more compact form. 
Note that there are
$14$ maximal chains in this interval.

\begin{figure}[htb]
$$
\epsfxsize=2.5in \epsfbox{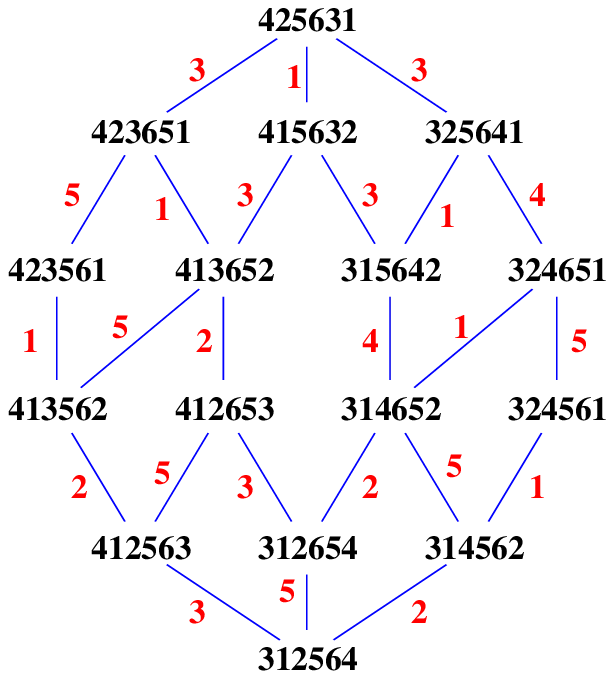}
$$
\end{figure}

\noindent Theorem~\ref{combiconst} gives us that
\begin{eqnarray*}
c_{(2,2,1)}^\zeta&=&\left|{\mathcal H}_{(2,2,1)}(\zeta)\right|
      - \left|{\mathcal H}_{(1,3,1)}(\zeta)\right|
      - \left|{\mathcal H}_{(2,0,3)}(\zeta)\right|
      + \left|{\mathcal H}_{(1,0,4)}(\zeta)\right|\\ 
      &&\qquad\qquad+ \left|{\mathcal H}_{(-1,3,3)}(\zeta)\right|
      - \left|{\mathcal H}_{(-1,2,4)}(\zeta)\right|.
\end{eqnarray*}
The sets ${\mathcal H}_{(-1,3,3)}(\zeta)$ and 
${\mathcal H}_{(-1,2,4)}(\zeta)$ are both empty since the
indices contains a negative component. 
Looking at Figure 2, we find 
$$
  {\mathcal H}_{(2,2,1)}(\zeta)=\{\hu{12}\hu{35}\hu{23}\hu{56}\hu{34},
  \hu{34}\hu{45}\hu{12}\hu{56}\hu{24}\}
$$ 
and ${\mathcal H}_{(1,3,1)}(\zeta)={\mathcal H}_{(2,0,3)}(\zeta)=\emptyset$. 
Hence $c_{(2,2,1)}^\zeta=2$. 
Now for $\lambda=(2,1,1,1)$ and $\sigma\in\Sym_4$, the sequences
$\lambda_\sigma$ that do not contains a negative component are 
$(2,1,1,1)$, $(2,1,0,2)$, $(2,0,2,1)$, $(2,0,0,3)$, $(0,3,1,1)$, 
$(0,3,0,2)$, $(0,0,4,1)$ and 
$(0,0,0,5)$. For our example, we have
$\left|{\mathcal H}_{(2,1,1,1)}(\zeta)\right|=5$,
$\left|{\mathcal H}_{(2,1,0,2)}(\zeta)\right|=2$,
$\left|{\mathcal H}_{(2,0,2,1)}(\zeta)\right|=2$ and all the others are empty.
Hence $c_{(2,1,1,1)}^\zeta=5-2-2=1$. Using (1.3) for this example, we get
$c_\lambda^\zeta=0$ for the other $\lambda$, since $14=5*2+4*1$ is the total
number of maximal chains.

With Theorem~\ref{combiconst}, we are able to show combinatorialy
many of the symmetries of the $c_\lambda^\zeta$ that were first shown 
using geometry in~\cite{BS97a}.
Let us start with symmetries derived from algebraic structures in the 
cohomology of the flag
manifolds. 
If $\zeta\in\Sym_n\subset\Sym_\infty$, let 
$\omega_0=(n,n-1,\ldots,3,2,1)$ be the longest element of
$\Sym_n$. As in Lemma~\ref{vsym}, the vertical and horizontal symmetries 
of (1.1) imply the following
lemmas.

\begin{lem}\label{SymmV} The following map (vertical symmetry) is a bijection
 \begin{eqnarray*}
  \Psi_v\colon\ R_{\u{}}(\zeta)\!\!\qquad\ \qquad&\longrightarrow&
    R_{\u{}}(\omega_0 \zeta\omega_0)\\ 
    \hu{\a_n\b_n}\cdots \hu{\a_1\b_1}\ &\longmapsto&\ 
    \hu{\omega_0(\b_n)\omega_0(\a_n)}\cdots \hu{\omega_0(\b_1)\omega_0(\a_1)}
 \end{eqnarray*}
\end{lem}

\begin{lem}\label{SymmH} 
The following map (horizontal symmetry) is a bijection
 \begin{eqnarray*}
 \Psi_h\colon\ R_{\u{}}(\zeta)\!\!\qquad\ \qquad&\longrightarrow&
    R_{\u{}}(\zeta^{-1})\\ 
    \hu{\a_n\b_n}\cdots \hu{\a_1\b_1}\ &\longmapsto&
    \hu{\a_1\b_1}\cdots \hu{\a_n\b_n}
 \end{eqnarray*}
\end{lem}
Given a non-degenerate hermitian form on $\C^n$, 
we get an involution on the flag manifold induced
by taking orthogonal complements. On the Schubert basis, 
this corresponds  to 
${\mathfrak S}_w\mapsto{\mathfrak S}_{\omega_0 w\omega_0}$ and 
${\mathfrak S}_{v(\lambda,k)}\mapsto{\mathfrak S}_{v(\lambda^t,n-k)}$, 
where $\lambda^t$ denotes the conjugate
partition. 
Thus $c_\lambda^\zeta=c_{\lambda^t}^{\omega_0\zeta\omega_0}$. 
Also $c_{u v(\lambda,k)}^w$ is the coefficient of 
${\mathfrak S}_{\omega_0}$ in ${\mathfrak
S}_{\omega_0w}{\mathfrak S}_u{\mathfrak S}_{v(\lambda,k)}$. 
Interchanging the roles of $u$ and $w$ gives
$c_\lambda^\zeta=c_{\lambda}^{\omega_0\zeta^{-1}\omega_0}$. 
Combining these, we get
$c_\lambda^\zeta=c_{\lambda^t}^{\zeta^{-1}}$.
Here we show these identities directly from Theorem~\ref{combiconst} 
and its dual version.

\begin{cor}\label{SymA} 
$c_\lambda^\zeta=c_{\lambda^t}^{\zeta^{-1}}$.
\end{cor}

\noindent{\em Proof  }
We note that both~(4.2) and the Jacobi identity have dual versions. 
For this, let $1^n$ denote the
partition conjugate to $(n)$, that is $(1,1,\ldots,1)$. 
From~\cite{Sottile96} or other formulations, the reader
deduces that
$$
 {\mathfrak S}_u{\mathfrak S}_{v(1^n,k)} = \sum_
      {\stackrel{\mbox{\scriptsize
 $\hx=\hu{\a_n\b_n}\cdots \hu{\a_1\b_1}\not\equiv 0$ }}
 {\a_1>\a_2>\cdots >\a_n}} {\mathfrak S}_{(\hx 1)u}. 
 \eqno{(4.5)}
$$
Here ${\mathfrak S}_{v(1^n,k)}=e_n(x_1,x_2,\ldots,x_k)$ is the 
$n$th elementary symmetric
polynomial. On the other hand, we know from~\cite{Macdonald91} that
$$ 
 {\mathfrak S}_{v(\lambda^t,k)}=S_{\lambda^t}(x_1,x_2,\ldots,x_k)=
      \det\left(e_{\lambda_i+j-i}(x_1,x_2,\ldots,x_k)
      \right)_{1\le i,j\le r}.
 \eqno{(4.6)}
$$
For $p\in\Z^r$, we define ${\mathcal E}_p(\zeta)=\emptyset$ if $p_i<0$ for
some $i$, and set  ${\mathcal E}_p(\zeta)$ to be 
$$
  \big\{\hx\in R_{\u{}}(\zeta) :
  \a_1>\cdots>\a_{p_1};\a_{p_1+1}>\cdots>\a_{p_1+p_2};
  \ldots;\a_{n-p_r+1}>\cdots>\a_{n}\big\}.
$$
otherwise. 
With computations similar to~(4.4), 
using a different expansion for the determinant~(4.6), we deduce that
$$  
 c_{\lambda^t}^\zeta=\sum_{\sigma\in\Sym_r}\epsilon(\sigma)
  \left|{\mathcal E}_{\overleftarrow{\lambda_\sigma}}(\zeta)\right|,
 \eqno{(4.7)}
$$
where for $p=(p_1,p_2,\ldots,p_r)$ we define 
$\overleftarrow{p}=(p_r,\ldots,p_2,p_1)$.
Now we note that $\Psi_h$ in Lemma~\ref{SymmH} maps 
${\mathcal H}_{p}(\zeta)$ bijectively to ${\mathcal
E}_{\overleftarrow{p}}(\zeta^{-1})$. 
Hence by Theorem~\ref{combiconst}, the equation~(4.7) is equal
to
$c_\lambda^\zeta$.\qed\medskip

\begin{cor}\label{SymB} 
$c_\lambda^\zeta=c_{\lambda^t}^{\omega_0\zeta\omega_0}$.
\end{cor}

\noindent{\em Proof  }
We only sketch the proof here since it is very similar to that of 
Corollary~\ref{SymA}. 
First, we use a different version of~(4.5), see~\cite{Sottile96}:
  $${\mathfrak S}_u{\mathfrak S}_{v(1^n,k)} = \sum_
      {\stackrel{\mbox{\scriptsize 
     $\hx=\hu{\a_n\b_n}\cdots \hu{\a_1\b_1}\not\equiv 0 $}}
      {\b_1>\b_2>\cdots >\b_n}}
        {\mathfrak S}_{(\hx 1)u}. 
$$
From this we define ${\mathcal E}'_p(\zeta)$ to be
$$
 \big\{\hx\in R_{\u{}}(\zeta) :
 \b_1>\cdots>\b_{p_1};\b_{p_1+1}>\cdots>\b_{p_1+p_2};
 \ldots;\b_{n-p_r+1}>\cdots>\b_{n}\big\}
$$
for $p\in\N^r$. 
We deduce that
$$
 c_{\lambda^t}^\zeta=\sum_{\sigma\in\Sym_r}
 \epsilon(\sigma)
 \left|{\mathcal E}'_{\lambda_\sigma}(\zeta)\right|.
 \eqno{(4.8)}
$$
Finally we note that $\Psi_v$ in Lemma~\ref{SymmV} maps 
${\mathcal H}_p(\zeta)$ bijectively to 
${\mathcal E}'_p(\zeta)$ and this concludes our proof.
\qed

For an integer $a$,
let $\phi_a\colon\N\to\N$ be defined by $\phi_a(i)=i$ if $i<a$, and
$\phi_a(i)=i+1$ if $i\ge a$. This map
$\phi_a$ induces an imbeding $\phi_a^*\colon\Sym_\infty\to\Sym_\infty$ where
$\zeta'=\phi_a^*(\zeta)$ is the unique permutation defined by $\zeta'(a)=a$ and
$\phi_a\circ\zeta'=\zeta\circ\phi_a$. If fact, $\phi_a$ also induces a
monomorphism 
$\phi_a^*\colon\Mon\to\Mon$ where $\phi_a^*(\hu{\a\b})=
\hu{\phi_a(\a)\phi_a(\b)}$. This is obvious
since the map $\phi_a^*$ sends generators to generators and preserves the
relations. 
This shows that $\phi_a^*\colon\Sym_\infty\to\Sym_\infty$ is also  
$\preceq$-order
preserving.

\begin{lem}\label{mapp} 
$\phi_a^*\colon R_{\u{}}(\zeta)\to R_{\u{}}(\phi_a^*(\zeta))$ is a bijection.
\end{lem}

As far as we know, the next corollary was first discovered in \cite{BS97a}
using geometry.

\begin{cor}[Theorem 5.1.1 of \cite{BS97a}]\label{SymC} 
$c_\lambda^\zeta=c_{\lambda}^{\phi_a^*(\zeta)}$.
\end{cor}

\noindent{\em Proof  }
This follows directly from Theorem~\ref{combiconst} since 
$\phi_a^*\big({\mathcal
H}_p(\zeta)\big)={\mathcal H}_p\big(\phi_a^*(\zeta)\big)$.
\qed\medskip

A.~Postnikov has communicated to us that he has also found 
combinatorial proofs of some of these
identities, particuliarly Proposition~\ref{XX} and Corollary~\ref{SymC}.

\section{Open Problems.}

One of the most enigmatic identities of \cite{BS97a} 
is the following proposition:

\begin{prop}[Theorem H of \cite{BS97a}]\label{cyclica} 
Let $\c=(1,2,3,\ldots,n)$ and $\zeta$ be in $\Sym_n$. 
Then
 $$
   c_{\lambda}^\zeta=c_{\lambda}^{\c \zeta \c^{-1}}.
 $$
\end{prop}

This was obtained using geometry, and as of now,
we do not know how to show this combinatorially. 
We note that (1.3) implies that 
$\big|{R_{\u{}}(\zeta)}\big|=\big|{R_{\u{}}(\c\zeta\c^{-1})}\big|$. 
This suggests the existence of a
bijection
$$
 \varphi\colon R_{\u{}}(\zeta)\longrightarrow R_{\u{}}(\c\zeta\c^{-1}).
 \eqno{(5.1)}
$$
Note that the two Posets
$[1,\zeta]_\preceq$ and $[1,\c\zeta\c^{-1}]_\preceq$ are not 
necessarily isomorphic.
 For example let $\zeta=(2,4,1,3)$, the interval 
$[1,\c\zeta\c^{-1}]_\preceq$ is a hexagon and
$[1,\zeta]_\preceq$ is not, it is a {\sl kite}.
For our next problem, we remark that the Jacobi identity~(4.3) 
is invertible, 
hence Proposition~\ref{cyclica} implies that 
$\big|{\mathcal H}_p(\zeta)\big|=\big|{\mathcal
H}_p(\c\zeta\c^{-1})\big|$ for any $p$.

\begin{prob}\label{Proa} 
Construct a bijection $\varphi$ as in~(5.1) 
such that 
$\varphi\big({\mathcal H}_p(\zeta)\big)={\mathcal H}_p(\c\zeta\c^{-1})$.
\end{prob}

A positive answer to this problem, combined with 
Theorem~\ref{combiconst}, would give a combinatorial
proof of Proposition~\ref{cyclica}.

Another direction of inquiry is to improve on Theorem~\ref{combiconst}.
It is a useful combinatorial description of the $c_\lambda^\zeta$ but it
is very unsatisfactory.
It would be more elegant to have a formula that does not involve signs.
Using symmetric polynomials~\cite{Macdonald91}, we know that
$$
	{\mathfrak S}_{v((p_1),k)}{\mathfrak S}_{v((p_1),k)}
	\cdots{\mathfrak S}_{v((p_1),k)} \ = \
	{\mathfrak S}_{v(\lambda,k)} +
	\sum_{\mu\triangleleft\lambda} a_{\lambda\mu} 
	{\mathfrak S}_{v(\mu,k)},
	\eqno{(5.2)}
$$
where $\lambda=\lambda(p)$ is the partition of $n$ obtained by 
rearranging the numbers
$p_1,p_2,\ldots,p_r$ in decreasing order, and
$\triangleleft$ is the strict dominance order on partitions.
Iterating~(4.2) and using~(5.2) it becomes clear
that $c_\lambda^\zeta\le\left|{\mathcal H}_p(\zeta)\right|$ 
for any $p$ such that $\lambda(p)=\lambda$. 
Now, let ${\mathcal D}_\lambda(\zeta)=\sum_{\sigma\in\Sym_r} 
{\mathcal H}_{\lambda_\sigma}(\zeta)$ where the
sum here denotes the disjoint union of sets. In general, 
a chain $\hx$ of $[1,\zeta]_\preceq$ weakly
fits many compositions of $n$. 
We use pairs $(\hx,\sigma)$ to describe elements in ${\mathcal
D}_\lambda(\zeta)$ meaning that
$\hx\in{\mathcal H}_{\lambda_\sigma}(\zeta)$. 
We extend the definition of $\epsilon$ to pairs
$(\hx,\sigma)$ by setting
$\epsilon(\hx,\sigma)=\epsilon(\sigma)$. 
The Equation~(4.4) can be rewritten as
  $$
c_\lambda^\zeta=\sum_{(\hx,\sigma)\in{\mathcal D}_\lambda} 
\epsilon(\hx,\sigma).
$$ 
In many cases, we can construct an involution 
$\theta\,\colon {\mathcal D}_\lambda(\zeta)\to{\mathcal D}_\lambda(\zeta)$ 
such that 
\begin{enumerate}
 \item[(i)] $\theta(\hx,\sigma)=(\hx,\sigma)$ only if $\sigma$ 
 is the identity, and
 \item[(ii)] if $\theta(\hx,\sigma)\neq(\hx,\sigma)$ then
 $\epsilon\big(\theta(\hx,\sigma)\big)=-\epsilon(\hx,\sigma)$.
\end{enumerate}

\noindent When this happens, we get a very nice combinatorial construction of 
$c_{\lambda}^\zeta$ since
$$
 c_\lambda^\zeta=\sum_{
 \stackrel{\mbox{\scriptsize$\hx\in{\mathcal H}_\lambda$}}%
 {\theta(\hx,1)=(\hx,1) }}
 1=\left|\{\hx\in{\mathcal H}_\lambda : \theta(\hx,1)=(\hx,1)\}\right|.
 \eqno{(5.3)}
$$

\begin{prob}\label{Prob} Find an involution 
$\theta\,\colon {\mathcal D}_\lambda(\zeta)\to{\mathcal
D}_\lambda(\zeta)$ for any $\zeta$ and $\lambda$.
\end{prob}

With $\lambda=(2,1,1,1)$ and $\zeta=(1,4,3,5,6,2)$, as in Figure 2, 
the set ${\mathcal
D}_\lambda(\zeta)$ contains nine elements. Here, it is clear how  to
construct $\theta$ 
and the only fixed point of $\theta$ is the chain
$\hu{12}\hu{34}\hu{45}\hu{56}\hu{24}$. 

It is interesting to note that among all the previously proposed 
conjectures to describe the $c_\lambda^\zeta$ combinatorially, none works. 
Using the monoid $\Mon$, it is
relatively easy to test them against Equation~(1.3). 
We use Proposition~\ref{Izeta} to get $u\le_k w=\zeta u$. 
We then use
Theorem~\ref{genere} to produce all maximal chains in
$[1,\zeta]_\preceq$ from the CM-chain. Finally we compare the two sides
of~(1.3). 
For example, in~\cite{Winkel_multiplication}, it is suggested that if we
display the numbers $\b_1$, $\b_2$, $\ldots$, $\b_n$ in a right adjusted
shape $\lambda$, then $c_\lambda^\zeta$ counts the number of chains that
weakly fit $\lambda$ and are strictly decreasing in every column (using the
$\b$'s). 
A small counterexample to this is obtained using $\zeta=(3,5,4,2,1)$. 

In the cases where we know how to construct the involution $\theta$, we have
used a Schensted-like insertion algorithm. 
This is an explicit correspondence ${\mathcal H}_{(1,n)}(\zeta)\longrightarrow
{\mathcal H}_{(n,1)}(\zeta)$. 
For this we consider the following transformations:
\begin{enumerate}
\item[{\bf A)}] 
$\hu{c\c}\hu{a\a}\hu{b\b}\mapsto\hu{a\a}\hu{c\c}\hu{b\b}$,\qquad 
         if $\{a,\a\}\cap\{c,\c\}=\emptyset$ and $a<b<c$,
\item[{\bf B)}] $\hu{\b\c}\hu{\a\b}\hu{\b\d}\mapsto\hu{\a\c}
\hu{\c\d}\hu{\b\c}$,\hskip1.65em
         if $\a<\b<\c<\d$,
\item[{\bf C)}] 
$\hu{\b\d}\hu{\a\b}\hu{\b\c}\mapsto\hu{\b\c}\hu{\c\d}\hu{\a\c}$,\hskip1.7em
         if $\a<\b<\c<\d$,
\item[{\bf D)}] 
$\hu{c\c}\hu{ac}\hu{b\b}\mapsto\hu{b\b}\hu{c\c}\hu{ac}$,\qquad \hskip.4em
         if $\{a,c,\c\}\cap\{b,\b\}=\emptyset$ and $a<b<c$,
\item[{\bf E)}] 
$\hu{b\b}\hu{ac}\hu{c\c}\mapsto\hu{ac}\hu{c\c}\hu{b\b}$,\hskip2.4em
         if $\{a,c,\c\}\cap\{b,\b\}=\emptyset$ and $a<b<c$,
\item[{\bf F)}] 
$\hu{b\b}\hu{a\a}\hu{c\c}\mapsto\hu{b\b}\hu{c\c}\hu{a\a}$,\hskip2.2em
         if $\{a,\a\}\cap\{c,\c\}=\emptyset$ and $a<b<c$,
\end{enumerate}

\smallskip
\noindent
The algorithm is very simple: To a chain in ${\mathcal H}_{(1,n)}(\zeta)$,
we keep on applying the 
transformations {\bf A} to {\bf F} to the rightmost triples. 
When we stop, we have a chain in ${\mathcal  H}_{(n,1)}(\zeta)$. 
The analysis of this algorithm
can be found at: \hfill\break 
{\tt http://www.math.yorku.ca/bergeron/appendix.html}.

There is one last identity of \cite{BS97a} that is asking for a
combinatorial proof. 
Given $\eta,\zeta\in\Sym_\infty$, we say that the permutations are
$\u{}$-disjoint if for any  
$\hu{\a_n\b_n}\cdots \hu{\a_1\b_1}\in R_{\u{}}(\zeta)$ and any
$\hu{\c_n\d_n}\cdots \hu{\c_1\d_1}\in 
R_{\u{}}(\eta)$ we have
$\{\a_1,\b_1,\ldots,\a_n,\b_n\}\cap\{\c_1,\d_1,\ldots,\c_n,\d_n\}=
\emptyset$ and
$\hu{\a_n\b_n}\cdots \hu{\a_1\b_1}\hu{\c_n\d_n}\cdots \hu{\c_1\d_1}\ne 0$.
From the relations (1.1) it is clear that this definition depends only on
one choice of an element 
from each $R_{\u{}}(\zeta)$ and $R_{\u{}}(\eta)$.

\begin{prop}[Theorem G of \cite{BS97a}]
\label{disjoint} 
Given $\eta$ and $\zeta$ $\u{}$-disjoint, we
have $$ c_{\lambda}^{\zeta\eta}=\sum_{\nu,\mu}
c_{\nu\mu}^{\lambda}c_{\nu}^{\zeta}c_{\mu}^{\eta},$$
where $c_{\nu\mu}^{\lambda}$ is the classical Littlewood-Richardson 
coefficient.
\end{prop}

Here it is not difficult to see how $R_{\u{}}(\zeta\eta)$ is related to 
$R_{\u{}}(\zeta)$ and
$R_{\u{}}(\eta)$. In fact $\u{}$-disjointness directly implies that
$[1,\zeta\eta]_\preceq\cong[1,\zeta]_\preceq\times [1,\eta]_\preceq$.
We can use that to relate ${\mathcal  H}_{p}(\zeta\eta)$ to 
${\mathcal  H}_{q}(\zeta)$ and ${\mathcal H}_{q'}(\eta)$ for some $q,q'$.

\begin{prob}\label{Proc} 
For $\zeta,\eta$ two $\u{}$-disjoint permutations, use
$[1,\zeta\eta]_\preceq\cong[1,\zeta]_\preceq\times [1,\eta]_\preceq$ and 
Theorem~\ref{combiconst} to
construct a combinatorial proof of Proposition~\ref{disjoint}.
\end{prob}

We end our list with problems related to $\Mon$ and $\preceq$.

\begin{prob}\label{Prod} 
Let $m=\ell_{\u{}}(\zeta)$. Equation~(1.3) suggests that we could:
\begin{enumerate}
\item[(a)] Find a representation of the symmetric group $\Sym_m$ on 
$\Q R_{\u{}}(\zeta)$ with
character given by $$\sum_{\lambda} c_\lambda^\zeta \chi^\lambda.$$
\item[(b)] Find a partition of  $R_{\u{}}(\zeta)$ similar to the one 
discussed after the relations~(2.2).
\end{enumerate}
\end{prob}

\begin{prob}\label{Proe} 
Describe the polynomial $\displaystyle P_n(t)=\sum_{\zeta\in\Sym_n}
t^{\ell_{\u{}}(\zeta)}$.
\end{prob}

\noindent
Here, let us list the first few of these polynomials:
$$
\begin{array}{l}
 P_1(t)=1,\qquad P_2(t)=1+t,\qquad P_3(t)=1+3t+2t^2,\\ 
 P_4(t)=1+6t+10t^2+6t^3+t^4,\\ 
 P_5(t)=1+10t+30t^2+40t^3+27t^4+10t^5+2t^6,\\ 
 P_6(t)=1+15t+70t^2+155t^3+195t^4+156t^5+86t^6+33t^7+ 8t^8 + t^9,\\ 
 P_7(t)=1+21t+140t^2+455t^3+875t^4+1120t^5+1038t^6+735t^7+ 406t^8\\  
         \qquad\qquad + 175t^9 +58t^{10}+14t^{11}+2t^{12},\\ 
 P_8(t)=1+28t+252t^2+1120t^3+2996t^4+5432t^5+7252t^6+7562t^7+ 6398t^8\\  
         \qquad\qquad  + 4492t^9 +2652t^{10}
         +1324t^{11}+556t^{12}+192t^{13}+52t^{14}+10t^{15}+t^{16}.
\end{array}
$$
It is also instructive to display the Poset $(\Sym_4,\preceq)$. 
Here we represent the permutations
using disjoint cycles notation.

\begin{figure}[htb]
$$
\epsfxsize=4.1in \epsfbox{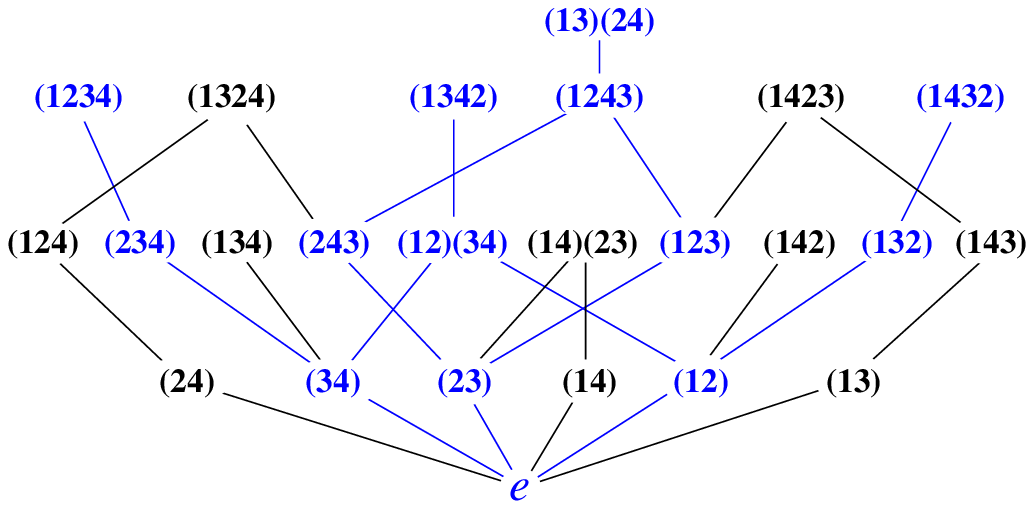}
$$
\end{figure}

\begin{prob}\label{Prof} 
What are the properties of the partial order $\preceq$. e.g. What is its
M\"obius function? Is any interval Cohen-Macauley? 
\end{prob}

We should mention here that the intervals contain {\sl hexagons} in general,
hence they are not 
shellable in the classical sense.

\begin{prob}\label{Prog} 
Is it possible to find a faithful representation of $\Mon$ as operators on
the polynomial ring $\Z[x_1,x_2,x_3,\ldots]$?
\end{prob}

This last problem is suggested by the situation for the nilplactic 
monoid $\mathcal N$. For $\mathcal N$ we
have  a faithful representation defined by $\u{i}\mapsto\partial_i$,  where
$\partial_i$ is the 
divided difference operator on $\Z[x_1,x_2,x_3,\ldots]$.

\bigskip
\noindent{\sc Acknowledgment} The authors are grateful to M. 
Shimozono and many others for
stimulating conversations.

\end{document}